\documentclass{article}

\usepackage{amssymb,latexsym,amsmath}

\usepackage{graphicx}

\begin{document}

\newcommand{\bfi}{\bfseries\itshape}

\makeatletter

\@addtoreset{figure}{section}

\def\thefigure{\thesection.\@arabic\c@figure}

\def\fps@figure{h, t}

\@addtoreset{table}{bsection}

\def\thetable{\thesection.\@arabic\c@table}

\def\fps@table{h, t}

\@addtoreset{equation}{section}

\def\theequation{\thesubsection.\arabic{equation}}

\makeatother

\newtheorem{thm}{Theorem}[section]

\newtheorem{prop}[thm]{Proposition}

\newtheorem{lema}[thm]{Lemma}

\newtheorem{cor}[thm]{Corollary}

\newtheorem{defi}[thm]{Definition}

\newtheorem{rk}[thm]{Remark}

\newtheorem{exempl}{Example}[section]

\newenvironment{exemplu}{\begin{exempl}  \em}{\hfill $\surd$

\end{exempl}}

\newcommand{\comment}[1]{\par\noindent{\raggedright\texttt{#1}

\par\marginpar{\textsc{Comment}}}}

\newcommand{\todo}[1]{\vspace{5 mm}\par \noindent \marginpar{\textsc{ToDo}}\framebox{\begin{minipage}[c]{0.95 \textwidth}

\tt #1 \end{minipage}}\vspace{5 mm}\par}

\newcommand{\ea}{\mbox{{\bf a}}}

\newcommand{\eu}{\mbox{{\bf u}}}

\newcommand{\ueu}{\underline{\eu}}

\newcommand{\ueo}{\overline{u}}

\newcommand{\oeu}{\overline{\eu}}

\newcommand{\ew}{\mbox{{\bf w}}}

\newcommand{\ef}{\mbox{{\bf f}}}

\newcommand{\eF}{\mbox{{\bf F}}}

\newcommand{\eC}{\mbox{{\bf C}}}

\newcommand{\en}{\mbox{{\bf n}}}

\newcommand{\eT}{\mbox{{\bf T}}}

\newcommand{\eL}{\mbox{{\bf L}}}

\newcommand{\eR}{\mbox{{\bf R}}}

\newcommand{\eV}{\mbox{{\bf V}}}

\newcommand{\eU}{\mbox{{\bf U}}}

\newcommand{\ev}{\mbox{{\bf v}}}

\newcommand{\eve}{\mbox{{\bf e}}}

\newcommand{\uev}{\underline{\ev}}

\newcommand{\eY}{\mbox{{\bf Y}}}

\newcommand{\eK}{\mbox{{\bf K}}}

\newcommand{\eP}{\mbox{{\bf P}}}

\newcommand{\eS}{\mbox{{\bf S}}}

\newcommand{\eJ}{\mbox{{\bf J}}}

\newcommand{\eB}{\mbox{{\bf B}}}

\newcommand{\eH}{\mbox{{\bf H}}}

\newcommand{\leb}{\mathcal{ L}^{n}}

\newcommand{\eI}{\mathcal{ I}}

\newcommand{\eE}{\mathcal{ E}}

\newcommand{\hen}{\mathcal{H}^{n-1}}

\newcommand{\eBV}{\mbox{{\bf BV}}}

\newcommand{\eA}{\mbox{{\bf A}}}

\newcommand{\eSBV}{\mbox{{\bf SBV}}}

\newcommand{\eBD}{\mbox{{\bf BD}}}

\newcommand{\eSBD}{\mbox{{\bf SBD}}}

\newcommand{\ecs}{\mbox{{\bf X}}}

\newcommand{\eg}{\mbox{{\bf g}}}

\newcommand{\paromega}{\partial \Omega}

\newcommand{\gau}{\Gamma_{u}}

\newcommand{\gaf}{\Gamma_{f}}

\newcommand{\sig}{{\bf \sigma}}

\newcommand{\gac}{\Gamma_{\mbox{{\bf c}}}}

\newcommand{\deu}{\dot{\eu}}

\newcommand{\dueu}{\underline{\deu}}

\newcommand{\dev}{\dot{\ev}}

\newcommand{\duev}{\underline{\dev}}

\newcommand{\weak}{\stackrel{w}{\approx}}

\newcommand{\mild}{\stackrel{m}{\approx}}

\newcommand{\strong}{\stackrel{s}{\approx}}

\newcommand{\weakdown}{\rightharpoondown}

\newcommand{\opg}{\stackrel{\mathfrak{g}}{\cdot}}

\newcommand{\opunu}{\stackrel{1}{\cdot}}
\newcommand{\opdoi}{\stackrel{2}{\cdot}}

\newcommand{\opn}{\stackrel{\mathfrak{n}}{\cdot}}

\newcommand{\tr}{\ \mbox{tr}}

\newcommand{\Ad}{\ \mbox{Ad}}

\newcommand{\ad}{\ \mbox{ad}}

\renewcommand{\contentsname}{ }

\title{Non maximal cyclically monotone graphs and construction of 
a bipotential for the Coulomb's dry friction law}

\author{Marius Buliga\footnote{"Simion Stoilow" Institute of Mathematics of the Romanian Academy,
 PO BOX 1-764,014700 Bucharest, Romania, e-mail: Marius.Buliga@imar.ro }, 
 G\'ery de Saxc\'e\footnote{Laboratoire de 
  M\'ecanique de Lille, UMR CNRS 8107, Universit\'e des Sciences et 
  Technologies de Lille,
 B\^atiment Boussinesq, Cit\'e Scientifique, 59655 Villeneuve d'Ascq cedex, 
 France, e-mail: gery.desaxce@univ-lille1.fr},
 Claude Vall\'ee\footnote{Laboratoire de M\'ecanique des Solides, UMR CNRS 6610, UFR SFA-SP2MI, Bd M. et P. Curie, t\'el\'eport 2, BP 30179, 86962 Futuroscope-Chasseneuil cedex,  France, e-mail: vallee@lms.univ-poitiers.fr}}

\date{This version: 10.01.2009}

\maketitle

{\bf MSC-class:} 49J53; 49J52; 26B25
\begin{abstract}
We show a surprising connexion between a property of the inf convolution of a
family of convex lsc functions  and the fact that  the 
intersection of maximal cyclically monotone graphs is the critical set of a
bipotential. 
 
We then  extend the results from \cite{bipo1} to bipotentials convex covers,
generalizing the notion of a bi-implicitly convex lagrangian cover.
 
As an application we prove that the bipotential related to Coulomb's friction
law  is related to a specific  bipotential  convex cover with the property that 
any graph of the cover is non maximal cyclically monotone. 
\end{abstract}

\section{Introduction}

$X$ and $Y$ are topological, locally convex, real vector spaces of dual 
variables $x \in X$ and $y \in Y$, with the duality product 
$\langle \cdot , \cdot \rangle : X \times Y \rightarrow \mathbb{R}$. 
We shall suppose that $X, Y$ have topologies compatible with the duality 
product, that is: any  continuous linear functional on $X$ (resp. $Y$) 
has the form $x \mapsto \langle x,y\rangle$, for some $y \in Y$ (resp. 
$y \mapsto \langle x,y\rangle$, for some  $x \in X$).
 
To  any convex and lsc function $\phi: X \rightarrow 
\bar{\mathbb{R}}$ we associate a function called 
 {\bf separable 
bipotential} $\displaystyle b: X \times Y \rightarrow 
\bar{\mathbb{R}}$, defined by the formula: 
$$ b(x,y) = \phi(x) + \phi^{*}(y) $$
(for the general notion of a bipotential see Definition \ref{def2}). 
Here the function $\displaystyle \phi^{*}: Y \rightarrow 
\bar{\mathbb{R}}$ is the Fenchel conjugate of $\phi$, 
defined by the expression (\ref{dfen}).

The function $b$ is obviously bi-convex and lsc in each
argument. By Fenchel inequality we have 
\begin{equation}
b(x,y) \geq \langle x, y \rangle
\label{genfen}
\end{equation}

Also the following string of equivalences is true: 
\begin{equation}
 y \in \partial b(\cdot,y) (x) \, \Longleftrightarrow \, x \in 
\partial b(x, \cdot)(y)  
\, \Longleftrightarrow b(x,y) = \langle x , y \rangle
\label{genequi}
\end{equation}
which is just another way of writing the well known string of equivalences
$$ y \in \partial \phi(x) \, \Longleftrightarrow \, x \in \partial \phi^{*}(y) 
\, \Longleftrightarrow \phi(x) + \phi^{*}(y) = \langle x , y \rangle$$


In mechanics subgradient inclusions  $y \in \partial \phi(x)$ are related to 
{\bf associated} constitutive laws \cite{Halp JM 75}. There exist also {\bf non-associated}  
constitutive laws which take the form $y \in \partial b(\cdot,y) (x)$ 
for functions $b$ of two variables, convex and lsc in each argument and 
satisfying (\ref{genfen}), 
(\ref{genequi}), which are called bipotentials. Such
an approach to non-associated constitutive laws has been first proposed in
\cite{saxfeng}. Examples of such  non associated 
constitutive laws are:
non-associated Dr\"ucker-Prager \cite{sax boush KIELCE 93}  and 
Cam-Clay models in soil mechanics \cite{sax BOSTON 95}, 
cyclic Plasticity (\cite{sax CRAS 92},\cite{bodo sax EJM 01}) 
and Viscoplasticity \cite{hjiaj bodo CRAS 00} of metals with non linear 
kinematical hardening rule, Lemaitre's damage law \cite{bodo}, the coaxial 
laws (\cite{sax boussh 2},\cite{vall leri CONST 05}). 

Of special interest to us is the formulation in terms of bipotentials of the 
 Coulomb's friction law \cite{saxfeng}, 
\cite{sax CRAS 92}, 
\cite{boush chaa IJMS 02}, \cite{feng hjiaj CM 06}, \cite{fort hjiaj CG 02}, 
\cite{hjiaj feng IJNME 04}, \cite{sax boush KIELCE 93}, \cite{sax feng IJMCM 98}, \cite{laborde}.

In  \cite{bipo1} we solved two key problems: (a) when the graph of a 
given multivalued operator can be expressed as the set of critical points of a 
bipotentials, and (b) a method of construction of a bipotential associated 
(in the sense of point (a)) to a multivalued, typically non monotone, operator. 
The main tool was the notion of {\bf convex lagrangian cover} of the graph of 
the multivalued operator, and a related notion of implicit convexity of this 
cover. 

The results of \cite{bipo1}  apply  only to bi-convex, bi-closed graphs (for
short BB-graphs) 
 admitting  at least one convex lagrangian cover 
by {\bf maximal cyclically monotone graphs}. 
This is a rather large class of graph of multivalued operators but  
important applications to the mechanics, such as the 
bipotential associated to contact with friction \cite{saxfeng}, are not in 
this class. 

This paper  is dedicated to the  extension of  the method 
presented in \cite{bipo1} to a more general class of BB-graphs. This is done 
in two steps. In the first step we prove Theorem \ref{thmsup}, the main result
of this paper.The result is that  the intersection of two maximal cyclically
monotone graphs  is the critical set of a strong 
bipotential if and only if a condition formulated in terms of 
the inf convolution of a family of convex lsc functions is true. 
In the second step we extend  the main result of \cite{bipo1} by replacing 
the notion of convex lagrangian cover with the one of a 
{\bf bipotential convex cover} (definition \ref{defcover}). In this way we are able 
to apply our results to the bipotential for the Coulomb's friction law.

\paragraph{Aknowledgements.} We express our thanks to the anonymous referee for 
detecting a flaw in a previous version of the paper, and for several valuable
questions and suggestions for improvement of this article. The first author  acknowledges partial support 
from the Romanian Ministry of Education and Research, through the grant 
CEX06-11-12/2006.

\section{Notations and Definitions}

For any convex and closed set $A \subset X$, its  indicator function,  
$\displaystyle \chi_{A}$, is defined by 
$$\chi_{A} (x) = \left\{ \begin{array}{ll}
0 & \mbox{ if } x \in A \\ 
+\infty & \mbox{ otherwise } 
\end{array} \right. $$
The indicator  function is convex and lsc. 

We use the notation: $\displaystyle \bar{\mathbb{R}} = \mathbb{R}\cup \left\{ +\infty \right\}$. 

Given a function 
$\displaystyle \phi: X \rightarrow  \bar{\mathbb{R}}$, 
the conjugate $\phi^{*}: Y \rightarrow \bar{\mathbb{R}}$ 
is defined by: 
\begin{equation}
\phi^{*}(y) = \sup \left\{ \langle y,x\rangle - 
\phi(x) \mid x \in X \right\} \ .
\label{dfen}
\end{equation}
The conjugate is always convex and lsc. 

We denote by $\Gamma(X)$ the class of convex and lsc 
functions $\displaystyle \phi: X \rightarrow \bar{\mathbb{R}}$. 
The  class of convex and lsc functions  
$\displaystyle \phi: X \rightarrow \bar{\mathbb{R}}$ with non-empty epigraph is denoted by 
$\displaystyle \Gamma_{0}(X)$. 

The effective domain of a function $\phi: X \rightarrow \bar{\mathbb{R}}$ is $dom \, \phi \ = \ \left\{ x \in X \mbox{ : } \phi(x) < + \infty\right\}$. 

For $X$ real locally convex topological vector space and $f,g: X
\rightarrow \bar{\mathbb{R}}$ 
the inf-convolution of $f$ and $g$ at $x \in X$ is defined by: 
\begin{equation}
f \, \square \, g \, (x) \ = \ \inf_{x_{1}+ x_{2} = x} \left[ f(x_{1}) + g(x_{2}) \right] 
\label{dinf}
\end{equation}

The subgradient of a function $\displaystyle \phi: X \rightarrow
\bar{\mathbb{R}}$ at a point $x \in X$ is the (possibly empty) set: 
$$\partial \phi(x) = \left\{ u \in Y \mid \forall z \in X  \  \langle z-x, u \rangle \leq \phi(z) - \phi(x) \right\} \  .$$ 
In a similar way is defined the subgradient of a function $\psi: Y \rightarrow \bar{\mathbb{R}}$ in a point $y \in Y$, as the set: 
$$\partial \psi(y) = \left\{ v \in X \mid \forall w \in Y  \  \langle v, w-y \rangle \leq \psi(w) - \psi(y) \right\} \ .$$

\begin{defi}   {\bf The graph of a multivalued operator} $\displaystyle A: X
\rightarrow 2^{Y}$ is the set 
$$Gr(A) \, = \, \left\{(x,y) \in X \times Y \mid y \in A(x) \right\}$$  
Any subset $M \subset X \times Y$ is the graph of a operator $\displaystyle A: X
\rightarrow 2^{Y}$. Associated to $M$ is  the  multivalued operator
$$\displaystyle  X \ni x  \mapsto m(x) \ = \ \left\{ y \in Y \mid (x,y) \in 
M \right\} \ .$$
The {\bf dual} operator is given by 
$$\displaystyle  Y \ni y  \mapsto m^{*}(y) \ = \ \left\{ x \in X \mid (x,y) \in M \right\} \ .$$
The {\bf domain} of the graph $M$ (or the effective domain of the associated
operator $m$) is the set 
$\displaystyle dom(M) = \left\{ x \in 
X \mid m(x) \not = \emptyset\right\}$. 

The {\bf image} of the graph $M$ is the set $\displaystyle im(M) = \left\{ 
y \in Y \mid m^{*}(y) \not = \emptyset\right\}$. 
\label{def1}
\end{defi}

For any $\phi \in \Gamma(X)$ we shall denote by $M(\phi)$ the graph: 
\begin{equation}
M(\phi) \ = \ \left\{ (x,y) \in X \times Y \ \mid \ \phi(x)+\phi^{*}(y) = \langle x, y \rangle \right\} \  .
\label{mphi}
\end{equation}
The operator associated to the graph $M(\phi)$ is $\partial \phi$. 
The dual operator associated to $M(\phi)$ is 
$\partial \phi^{*}$ (the subgradient of the Legendre-Fenchel dual of $\phi$).

\begin{defi} A graph $M$ is {\bf cyclically monotone} if for all integer $m>0$
and any family of couples $(x_{j},y_{j}) \in M, j=0,1,\ldots,m$, 
\begin{equation}
\displaystyle  \langle x_{0} -  x_{m}, y_{m} \rangle + \sum_{k=1}^{m} \langle x_{k} -  x_{k-1}, y_{k-1} \rangle \leq 0.
\label{Cyclically Monotone inequality}
\end{equation}
A cyclically monotone graph $M$ is {\bf maximal} if it does not admit a strict prolongation which is cyclically monotone.
\label{def CM graph}
\end{defi}
By reindexing the couples, we easily recast the previous inequality as
\begin{equation}
\langle x_{m}, y_{0} - y_{m} \rangle
+ \sum_{k=1}^{m} \langle x_{k-1}, y_{k} - y_{k-1} \rangle \leq 0,
\label{Cyclically Monotone inequality 2}
\end{equation}
fact which shows that the graphs of a law and of its dual law are simultaneously cyclically monotone.
Rockafellar \cite{rocka} Theorem 24.8 (see also Moreau \cite{moreau} 
Proposition 12.2) proved a Theorem that can be stated as: 
\begin{thm}
Given a graph $M$, there exist a potential $\phi\in\Gamma_0(X)$ such that 
$M \subset Gr(\partial \phi)$ if and only if $M$ is cyclically monotone. 
The potential $\phi$ is unique up to an additive constant and it is 
 defined by
\begin{equation}
\displaystyle  \phi (x) =\ \sup \left\{   \langle x -  x_{m}, y_{m} \rangle + 
\sum_{k=1}^{m} \langle x_{k} -  x_{k-1}, y_{k-1} \rangle \right\}  + \phi (x_{0}  ),
 \label{phi Rockafellar}
\end{equation}
where $x_{0}$ and $\phi (x_{0})$ are arbitrarily fixed and the 'sup' is extended to any  $m>0$ and to any couples $(x_{k},y_{k}) \in M, k=1,2,\ldots,m$.
\label{thm Rockafellar}
 \end{thm}


Because the dual law is also cyclically monotone, we can apply once again  the construction of the previous Theorem, giving the function 
\begin{equation}
\psi (y) =\ \sup \left\lbrace  \langle x_{m}, y - y_{m} \rangle + \sum_{k=1}^{m} \langle x_{k-1}, y_{k} - y_{k-1} \rangle
\right\rbrace  + \psi (y_{0}  ),
 \label{psi Rockafellar}
\end{equation}
such that  $M \subset M(\psi^{*})$. Excepted when $M$ is maximal, $\phi$ and $\psi^{*}$ are in general distinct function, as it will be seen further in the application.


\begin{defi} A {\bf bipotential} is a function $b: X \times Y \rightarrow
 \bar{\mathbb{R}}$, with the properties: 
\begin{enumerate}
\item[(a)] $b$ is convex and lower semicontinuos in each argument; 
\item[(b)] for any $x \in X , y\in Y$ we have $\displaystyle b(x,y) \geq \langle x, y \rangle$; 
\item[(c)]  for any $(x,y) \in X \times Y$ we have the equivalences: 
\begin{equation}
y \in \partial b(\cdot , y)(x) \ \Longleftrightarrow \ x \in \partial b(x, \cdot)(y)  \ \Longleftrightarrow \ b(x,y) = 
\langle x , y \rangle \ .
\label{equiva}
\end{equation}
\end{enumerate}
The {\bf graph} of $b$ is 
\begin{equation}
M(b) \ = \ \left\{ (x,y) \in X \times Y \ \mid \ b(x,y) = \langle x, y \rangle \right\} \  .
\label{mb}
\end{equation}
\label{def2}
\end{defi}

Particular cases of bipotentials are separable ones, described in the
introduction. Many other non separable bipotentials exist. 

We introduce next the notion of a strong bipotential. Conditions (B1S) and (B2S) 
appear  as relations (51), (52) in \cite{laborde}. 

\begin{defi}
 A function $b: X \times Y \rightarrow 
\bar{\mathbb{R}}$ is a {\bf strong bipotential} if it satisfies the conditions: 
\begin{enumerate}
\item[(a)] $b$ is convex and lower semicontinuos in each argument; 
\item[(B1S)] for any  $y \in Y$ $\inf \left\{  b(z,y) - \langle z, y \rangle \mbox{ : } z \in X \right\} \in  \left\{ 0, + \infty \right\}$; 
\item[(B2S)] for any  $x \in X$ $\inf \left\{  b(x,p) - \langle x, p \rangle \mbox{ : } p \in Y \right\} \in  \left\{ 0, + \infty \right\}$. 
\end{enumerate}
\label{defstrongb}
\end{defi}

\begin{prop}
Any strong bipotential is a bipotential. 
\end{prop}

\paragraph{Proof.} 
Let $b$ be a strong bipotential. From (B1S) and (B2S) we have to prove conditions (b) and (c) of Definition \ref{def2}. 
  
Remark first that any of the two conditions (B1S) and (B2S) implies (b). All is left to prove is (c).

For this  take $x \in X$, $y \in Y$ such that  $y \in \partial 
b(\cdot , y)(x)$. This is equivalent with: $x$ is a minimizer of the function 
$\displaystyle z \in X \mapsto \left( b(z,y) - \langle z, y \rangle \right)$. But according 
to (B1S) the minimum of this function is equal to $0$. Therefore 
$b(x,y) = \langle x , y \rangle$. We proved that $y \in \partial 
b(\cdot , y)(x) \ 
\Longrightarrow \ b(x,y) =  \langle x , y \rangle$. The inverse implication is 
trivial, thus we have an equivalence. 

In the same, using (B2S)  we prove that $x \in \partial b(x , \cdot)(y)  \ \Longleftrightarrow \ b(x,y) = 
\langle x , y \rangle$ . The condition (c) is therefore satisfied and $b$ is a bipotential. 
\quad $\blacksquare$

Any separable bipotential is a strong bipotential. Indeed, for any convex lsc 
$\phi: X \rightarrow \bar{\mathbb{R}}$ we have: 
$$\inf \left\{  \phi(z) + \phi^{*}(y) - \langle z, y \rangle \mbox{ : } z \in X
\right\} \, \in \, \left\{ 0, + \infty\right\}$$ 
by the definition (\ref{dfen}) of $\displaystyle \phi^{*}$. The notion of a
strong bipotential (introduced in relations (51), (52)  \cite{laborde}) is
motivated also by the fact that all bipotentials considered in applications in
mechanics are in fact strong bipotentials.

Not all bipotentials are strong bipotentials. Consider for example $X = Y = \mathbb{R}$, with 
the duality $\langle x, y \rangle = xy$, and $b: \mathbb{R} \times \mathbb{R} 
\rightarrow \mathbb{R}$ defined by: 
$$b(x,y) \ = \ \mid x \mid \, \left( e^{-y} + 1 \right) + xy$$ 
This is bipotential which is not a strong bipotential, and $M(b) = \left\{ 0 \right\} \times
\mathbb{R}$. Indeed, $b$ is convex and lsc in each argument and 
$b(x,y) \geq xy$ for any $x, y \in \mathbb{R}$. It is easy to check that 
$y \in \partial b(\cdot , y)(x) \ \Longleftrightarrow \ x \in \partial b(x, \cdot)(y)  \ 
\Longleftrightarrow \ x \ = 0$. But $x = 0$ is equivalent with $b(x,y) = xy$, therefore 
$b$ is a bipotential. Nevertheless, this is not a strong bipotential. Indeed, for 
$x \not = 0$ we have 
$$\inf \left\{  b(x,p) - \langle x, p \rangle \mbox{ : } p \in \mathbb{R} \right\} = 
\mid x \mid \ \not \in \ \left\{ 0, + \infty \right\} $$

\section{Bipotentials for cyclically monotone graphs}

Maximal cyclically monotone graphs are critical sets of separable bipotentials.


The following Theorem shows that there exist bipotentials $b$ with the 
property that $M(b)$ is a {\bf cyclically monotone, but not maximal}  
set. In this section we show a surprising connection between 
bipotentials and the inf convolution operation.

\begin{thm}
Let  $b_{1}$ and $b_{2}$ be separable bipotentials associated respectively to 
the convex and lsc functions 
$\displaystyle \phi_{1} , \phi_{2}: X
\rightarrow \bar{\mathbb{R}}$, that is 
$$b_{i} (x,y) \ = \ \phi_{i}(x) + \phi_{i}^{*}(y)$$ 
for any $i = 1,2$ and $(x,y) \in X \times Y$. Consider the following 
assertions:
\begin{enumerate}
\item[(i)] 
$b= max (b_{1},b_{2})$ is a strong 
bipotential.
\item[(ii')] For any $\displaystyle y \in \, dom \,  \phi_{1}^{*} \, \cap \, dom \,  \phi_{2}^{*}$ and for any $\lambda \in 
[0,1]$ we have 
\begin{equation}
\left( \lambda \, \phi_{1} \, + \, (1-\lambda) \, \phi_{2} \right)^{*}(y) \ = 
\ \lambda \, \phi_{1}^{*}(y) \, + \, (1-\lambda) \, \phi_{2}^{*}(y) 
\label{iiprim}
\end{equation}
\item[(ii")] For any $\displaystyle x \in \, dom \,  \phi_{1} \, \cap \, dom \,  \phi_{2}$ and for any $\lambda \in 
[0,1]$ we have 
\begin{equation}
\left( \lambda \, \phi_{1}^{*} \, + \, (1-\lambda) \, \phi_{2}^{*} 
\right)^{*}(x) \ = 
\ \lambda \, \phi_{1}(x) \, + \, (1-\lambda) \, \phi_{2}(x)
\label{iisec}
\end{equation}
\end{enumerate} 
Then the point (i) is equivalent with the conjunction of (ii'), (ii"), (for
short: 
  (i) $\Longleftrightarrow$ ( (ii') AND (ii") ) ). 
\label{thmsup}
\end{thm}

\begin{rk}
If $\displaystyle b_{1}, b_{2}$ are separable bipotentials and 
$\displaystyle b = max (b_{1}, b_{2})$ is a bipotential then 
$\displaystyle M(b) = M(b_{1}) \cap M(b_{2})$, therefore $M(b)$ is 
the intersection of two maximal cyclically monotone graphs. 
\end{rk}

\paragraph{Proof.} 

Before we begin to prove the equivalence  let us remark that for any 
$x \in X$ and for any $y \in Y$,  the functions $b(\cdot, x)$ and $b(\cdot, y)$ are 
convex and 
lsc as 
superior envelopes of such functions. Also,  for any $(x,y) \in X \times Y$ 
we have $\displaystyle b_{1}(x,y)\geq \langle x,y \rangle$ and $\displaystyle 
b_{2}(x,y)\geq \langle x,y \rangle$, therefore 
$$\displaystyle b (x,y) = max \left( b_{1}(x,y), b_{2} (x,y) \right) \geq 
\langle x,y \rangle \quad .$$ 

Let $M(b)\subset X \times Y$ be the set of pairs $(x,y)$ with the property 
$b(x,y) = \langle x , y \rangle$. 
If $(x,y) \in M(b)$, then
$$ \langle x,y \rangle \leq b_{i} (x,y) \leq b (x,y) = \langle x,y 
\rangle \quad   (i=1,2)$$
which proves that $\displaystyle (x,y) \in M(b_{1})\cap M(b_{2})$. 
Conversely, if 
$\displaystyle (x,y) \in M(b_{1})\cap M(b_{2})$ then 
$\displaystyle   \langle x,y \rangle = b_{1} (x,y) = b_{2} (x,y) = b(x,y)$,
therefore  $(x,y) \in M(b)$. In conclusion $\displaystyle M(b) = M(b_{1}) \cap 
M(b_{2})$.

Thus the equivalence we have to prove becomes: 
\begin{enumerate}
\item[(I)] the condition (B1S) from Definition \ref{defstrongb}  is equivalent with (ii'), 
 \item[(II)]  the condition (B2S) from Definition \ref{defstrongb}  is equivalent with (ii"). 
\end{enumerate}

These two equivalences have similar proofs. We shall give the proof of the 
first equivalence. 

The function $b$ admits the following characterization: 
$$b(x,y) = \max_{\lambda \in [0,1]} \left\{ \lambda b_{1}(x,y) + (1-\lambda) 
b_{2}(x,y) \right\}$$
For $\lambda \in [0,1]$ denote by $\displaystyle b^{\lambda}(x,y) = 
\lambda b_{1}(x,y) + (1-\lambda) 
b_{2}(x,y)$. For any $y \in Y$ such that $\displaystyle \phi_{1}^{*}(y) < +
\infty$, $\displaystyle \phi_{2}^{*}(y) < +\infty$ define the set 
$$C(y) = \left\{ z \in X \mbox{ : } b_{1}(z,y) < + \infty \, , \,  
b_{2}(z,y) < + \infty \right\}$$ 
and remark that $C(y) \subset X$ is a convex set. In fact $\displaystyle 
C(y) = \, dom \, \phi_{1} \, \cap \, dom \, \phi_{2}$, therefore we may drop 
the $y$ argument and write $C$ instead of $C(y)$. 


Consider then the 
function $f(\cdot , \cdot , y): C \times [0,1] \rightarrow \mathbb{R}$ given
by $\displaystyle f(z,\lambda, y) = \langle z, y \rangle - b^{\lambda}(z,y)$. 
This function is affine and continuous in $\lambda$, $[0,1]$ is a compact 
convex subset of the vector space $\mathbb{R}$. Also, this function is concave 
and upper semicontinuous in $z \in C$. Therefore we are in position to 
apply the minimax Theorem of Sion \cite{sion} and deduce that: 


\begin{equation}
\min_{\lambda \in [0,1]} \sup_{z \in C} f(z, \lambda, y) \ = \ 
\sup_{z \in C} \min_{\lambda \in [0,1]} f(z, \lambda, y)
\label{minimax}
\end{equation}
Let us compute the terms of the equality (\ref{minimax}). We have: 
$$ A = \min_{\lambda \in [0,1]} \sup_{z \in C} f(z, \lambda, y) \ = \ 
\min_{\lambda \in [0,1]} \sup_{z \in C} \left\{ \langle z, y \rangle -
b^{\lambda}(z,y) \right\} \ = \ $$ 
$$ = \ \min_{\lambda \in [0,1]} \ \left( \lambda \, \phi_{1} \, + \, (1-\lambda) \, \phi_{2} 
\right)^{*}(y) \ - 
\ \lambda \, \phi_{1}^{*}(y) \, - \, (1-\lambda) \, \phi_{2}^{*}(y)$$
For the other term of the equality (\ref{minimax}) we have: 
$$B = \sup_{z \in C} \min_{\lambda \in [0,1]} f(z, \lambda, y) \ = \ 
\sup_{z \in C} \left\{ \langle z, y\rangle - b(z,y) \right\}$$
We have $A = B$ thus (\ref{minimax}) is equivalent with: 
\begin{equation}
\sup_{z \in C} \left\{ \langle z, y\rangle - b(z,y) \right\} \ = \ \min_{\lambda \in [0,1]} \ 
\left( \lambda \, \phi_{1} \, + \, (1-\lambda) \, \phi_{2} 
\right)^{*}(y) \ - 
\ \lambda \, \phi_{1}^{*}(y) \, - \, (1-\lambda) \, \phi_{2}^{*}(y) 
\label{equivalence}
\end{equation}

Suppose that $b$ is a strong bipotential and  let $y \in Y$ such that $\displaystyle \phi_{1}^{*}(y) < +
\infty$, $\displaystyle \phi_{2}^{*}(y) < +\infty$. This implies, by (B1S), that 
$$\sup_{z \in C} \left\{ \langle z, y\rangle - b(z,y) \right\} \ = \ 0$$ 
By (\ref{equivalence}) we deduce that 
$$\min_{\lambda \in [0,1]} \ 
\left( \lambda \, \phi_{1} \, + \, (1-\lambda) \, \phi_{2} 
\right)^{*}(y) \ - 
\ \lambda \, \phi_{1}^{*}(y) \, - \, (1-\lambda) \, \phi_{2}^{*}(y)  \ = \ 0$$
But in general we have 
$$\max_{\lambda \in [0,1]} \ 
\left( \lambda \, \phi_{1} \, + \, (1-\lambda) \, \phi_{2} 
\right)^{*}(y) \ - 
\ \lambda \, \phi_{1}^{*}(y) \, - \, (1-\lambda) \, \phi_{2}^{*}(y) \ \leq \ 0$$
which comes from the equality true for any $\lambda \in [0,1]$: 
$$\langle z, y\rangle - \lambda \, \phi_{1}(z) - \, (1-\lambda) \, \phi_{2}(z) \
 = \ \lambda \, \left( \langle z, y \rangle - \phi_{1}(z) \right) \, + \, (1-\lambda) \, \left(  \langle z, y \rangle -
 \phi_{2}(z) \right)$$ 
We find therefore that for any $\lambda \in [0,1]$ 
$$\left( \lambda \, \phi_{1} \, + \, (1-\lambda) \, \phi_{2} 
\right)^{*}(y) \ - 
\ \lambda \, \phi_{1}^{*}(y) \, - \, (1-\lambda) \, \phi_{2}^{*}(y)  \ = 0$$

Conversely, let $y \in Y$. If $\displaystyle \phi_{1}^{*}(y) = +
\infty$ or  $\displaystyle \phi_{2}^{*}(y) = +\infty$ then we have: 
$$\inf \left\{  b(z,y) - \langle z, y \rangle \mbox{ : } z \in X \right\} \ = \ + \infty$$
Suppose that $\displaystyle \phi_{1}^{*}(y) < +
\infty$ and  $\displaystyle \phi_{2}^{*}(y) < +\infty$. From (ii') and (\ref{equivalence}) we deduce that 
$$\inf \left\{  b(z,y) - \langle z, y \rangle \mbox{ : } z \in X \right\} \ = \ 0$$ 
The second equivalence has a similar proof. \quad $\blacksquare$

It is easy to construct examples of separable bipotentials $\displaystyle b_{i} (x,y) \ = \ \phi_{i}(x) + 
\phi_{i}^{*}(y)$,  $i = 1,2$,  such that $\displaystyle b = \max \left\{b_{1}, b_{2} \right\}$ is not a bipotential. 
For this let us take $X = Y = \mathbb{R}$ with the duality given by the product, and let us choose 
$\displaystyle \phi_{1}, \phi_{2}$ smooth (for example $\displaystyle \mathcal{C}^{1}$) convex functions defined on 
$\mathbb{R}$ with values in $\mathbb{R}$. Then $\displaystyle M(b_{1}) =
Gr(\phi_{1}')$  and $\displaystyle M(b_{2}) =
Gr(\phi_{2}')$, therefore $M(b)$ is the intersection of the graphs of 
$\displaystyle \phi_{1}'$ and $\phi_{2}'$. In general  $M(b)$ is  
not bi-convex, because it is just an intersection of graphs of increasing continuous functions. For example, if 
$\displaystyle \phi_{1}(x) = x^{4}/4$ and $\displaystyle \phi_{2}(x) = x^{2}/2$ then 
$M(b) = \left\{ (1,1) , (0,0), (-1, -1) 
\right\}$, which is not bi-convex.

The conditions (ii'), (ii") from Theorem \ref{thmsup} imply relations which 
can be expressed with the
help of inf convolutions. Let us examine condition (ii'); the same arguments can
be used for the symmetric condition (ii'). Consider $\displaystyle \phi_{1},
\phi_{2} \in \Gamma_{0}(X)$. For any $\lambda \in (0,1)$ we introduce two
functions defined on $X$ by: 
$$f_{1, \lambda} (x) \, = \, \lambda \, \phi_{1}(\frac{1}{\lambda} x) \quad , 
\quad f_{2, \lambda} (x) \, = \, (1-\lambda) \, \phi_{2}(\frac{1}{1-\lambda} x)$$
Then condition (ii") implies that for any $\displaystyle x \in \, dom \,  \phi_{1} \, \cap \, dom \,  \phi_{2}$ and for any $\lambda \in 
(0,1)$ we have 
\begin{equation}
f_{1, \lambda}(\lambda x) \, + \, f_{2,\lambda} ((1-\lambda) x)  \ = 
\  f_{1,\lambda} \square f_{2,\lambda} (x)
\label{iisecconv}
\end{equation}
Indeed, we have 
$$\left(f_{1,\lambda}^{*} + f_{2,\lambda}^{*}\right)^{*}(x) \, = \, 
\left( \lambda \, \phi_{1}^{*} \, + \, (1-\lambda) \, \phi_{2}^{*} 
\right)^{*}(x)$$
The space $X$ is locally convex therefore 
$$\left(f_{1,\lambda}^{*} + f_{2,\lambda}^{*}\right)^{*}(x) \, \leq \, 
 f_{1,\lambda} \square f_{2,\lambda} (x)$$
Therefore   (ii") implies that
 $$f_{1, \lambda}(\lambda x) \, + \, f_{2,\lambda} ((1-\lambda) x)  \, \leq  
\,  f_{1,\lambda} \square f_{2,\lambda} (x)$$
which, by the definition (\ref{dinf}) of the inf convolution operation, is
equivalent with (\ref{iisecconv}).

This is leading us to the following corollary of Theorem \ref{thmsup}. 
\begin{cor}
Let $\displaystyle \phi_{1},
\phi_{2} \in \Gamma_{0}(X)$ such that 
$$b(x,y) \, = \, \max \left(\phi_{1}(x) + \phi_{1}^{*}(y), \phi_{2}(x) +
\phi_{2}^{*}(y) \right\}$$ 
is a strong bipotential. Then, with the previous notations, for any  
$\displaystyle x \in \, dom \,  \phi_{1} \, \cap \, dom \,  \phi_{2}$ and for any $\lambda \in 
(0,1)$ we have 
$$\partial \left( f_{1,\lambda} \square f_{2,\lambda} \right) (x) \, = \, 
\partial \phi_{1} (x) \, \cap \, \partial \phi_{2} (x)$$
\label{pconvo}
\end{cor}

\paragraph{Proof.}
By Theorem \ref{thmsup}, if $b$ is a strong bipotential then (ii") is true. 
By previous reasoning this implies the relation (\ref{iisecconv}). We apply then
 Lemma 2.6, Lemma 2.7 \cite{zagrod} and we obtain that  
$$\partial \left( f_{1, \lambda} \, \square \,  f_{2, \lambda} \right) (x) \ = \ \partial
f_{1,\lambda}(\lambda x) \, \cap \, \partial f_{2, \lambda} ((1-\lambda) x)$$
From the definition of the functions $\displaystyle f_{1, \lambda}, f_{2,
\lambda}$ we obtain the conclusion.  $\blacksquare$

The  following interesting  question has been suggested by the anonymous 
referee:  can the maximum of two separable bipotentials  be
a non strong bipotential?

\section{Bipotential convex covers}

Let $Bp(X,Y)$ be the set of all bipotentials $b: X \times Y \rightarrow
\bar{\mathbb{R}}$. We shall need the following Definition concerning implicitly convex functions.

\begin{defi}
Let $\Lambda$ be an arbitrary non empty set and $V$ a real vector space. The 
function $f:\Lambda\times V \rightarrow \bar{\mathbb{R}}$ is 
{\bf implicitly  convex} if for any two elements 
$\displaystyle (\lambda_{1}, z_{1}) , 
(\lambda_{2},  z_{2}) \in \Lambda \times V$ and for any two numbers 
$\alpha, \beta \in [0,1]$ with $\alpha + \beta = 1$ there exists 
$\lambda  \in \Lambda$ such that 
$$f(\lambda, \alpha z_{1} + \beta z_{2}) \ \leq \ \alpha 
f(\lambda_{1}, z_{1}) + \beta f(\lambda_{2}, z_{2}) \quad .$$
\label{defimpl}
\end{defi}

In the following Definition we generalize the notion of a {\bf bi-implicitly 
convex lagrangian cover}, Definitions 4.1 and 6.6 \cite{bipo1}.

\begin{defi}   A {\bf bipotential
convex cover} of the non empty set $M$ is a function   
$\displaystyle \lambda \in \Lambda \mapsto b_{\lambda}$ from  $\Lambda$ with 
values in the set  $Bp(X,Y)$, with the 
properties:
\begin{enumerate}
\item[(a)] The set $\Lambda$ is a non empty compact topological space, 
\item[(b)] Let $f: \Lambda \times X \times Y \rightarrow \mathbb{R} \cup
\left\{ + \infty \right\}$ be the function defined by 


$$f(\lambda, x, y) \ = \ b_{\lambda}(x,y) .$$


Then for any $x \in X$ and for any $y \in Y$ the functions 
$f(\cdot, x, \cdot): \Lambda \times Y \rightarrow \bar{\mathbb{R}}$ and 
$f(\cdot, \cdot , y): \Lambda \times X \rightarrow \bar{\mathbb{R}}$ are  lower 
semi continuous  on the product spaces   $\Lambda \times Y$ and respectively 
$\Lambda \times X$ endowed with the standard topology, 
\item[(c)] We have $\displaystyle M  \ = \  \bigcup_{\lambda \in \Lambda} 
M(b_{\lambda})$.
\item[(d)] with the notations from point (b), the functions $f(\cdot, x, \cdot)$ 
and $f(\cdot, \cdot , y)$ are implicitly convex in the sense of Definition 
\ref{defimpl}.
\end{enumerate}
\label{defcover}
\end{defi}

Several remarks are in order. 

\begin{rk} A bipotential convex cover  $\displaystyle \lambda \in \Lambda 
\mapsto b_{\lambda}$ such that for any $\lambda \in \Lambda$ the bipotential 
$\displaystyle b_{\lambda}$ is separable is a bi-implicitly 
convex lagrangian cover. For such covers the sets $\displaystyle M(b_{\lambda})$
are {\bf maximal cyclically monotone} for any $\lambda \in \Lambda$.
\end{rk}

\begin{rk}
In general bipotential convex covers are {\bf not lagrangian} (see remark 6.1 
\cite{bipo1} for a justification of the "lagrangian" term). In the language of 
convex analysis this means that the sets $\displaystyle M(b_{\lambda})$ are 
not supposed to be cyclically monotone. 
\end{rk}

We shall see in the section concerning the applications to the Coulomb's friction
law that there exists bipotential convex covers with the property that for any 
$\lambda \in \Lambda$ the set $\displaystyle M(b_{\lambda})$ is {\bf cyclically
monotone but non maximal}. This is done by using bipotential covers constructed 
with the help of Theorem \ref{thmsup}. 

A bipotential convex cover is in some sense described by the collection 
$\displaystyle \left\{ b_{\lambda} \mbox{ : } \lambda \in \Lambda \right\}$. 
This is shown in the next Proposition. 


\begin{prop}
Let $\displaystyle \lambda \in \Lambda \mapsto b_{\lambda} \in Bp(X,Y)$ be a
bipotential convex cover and $g: \Lambda \rightarrow \Lambda$ be a continuous, 
invertible, with continuous inverse, function. Then 
 $\displaystyle \lambda \in \Lambda \mapsto b_{g(\lambda)} \in Bp(X,Y)$ is a
 bipotential convex cover. 
\end{prop}


\paragraph{Proof.}
This is obvious due to the general fact  that if 
$f:\Lambda\times V \rightarrow \bar{\mathbb{R}}$ is implicitly convex and 
$g: \Lambda \rightarrow \Lambda$ is a bijection then 
the function $f':\Lambda\times V \rightarrow \bar{\mathbb{R}}$, 
$\displaystyle f'(\lambda, x) = f(g(\lambda), x)$ is implicitly convex. 
This reflects into the fact that a bipotential convex cover is a notion 
 invariant with respect to {\bf continuous}
reparametrizations of $\Lambda$ (the continuity is needed in order to preserve 
the lower semi continuity assumptions from point (b) of the Definition
\ref{defcover}).  \quad $\blacksquare$ 

The next Theorem generalizes Theorem 6.7, the main result of \cite{bipo1}. We
shall skip its proof because it is just a rephrasing of the proof of Theorem 
6.7 \cite{bipo1}.

\begin{thm} Let $\displaystyle \lambda \mapsto b_{\lambda}$ be a bipotential 
convex cover  of 
 the graph $M$ and $b: X \times Y \rightarrow R$ defined by
\begin{equation}
b(x,y) \ = \ \inf \left\{ b_{\lambda}(x,y) \ \mid \  \lambda \in \Lambda \right\} \ . 
\end{equation}
Then $b$ is a bipotential and $M=M(b)$. 
\label{thm2}
\end{thm}

\section{Application: Coulomb's law of dry friction contact}

 This is is a typical 
example of what is called a non associated constitutive law in mechanics. 
Despite of its rather complex structure, it is worthwhile to have interest 
in it because of its importance in many practical problems. 

We shall not discuss here the phenomenal and experimental aspects but only 
the mathematical modeling with respect to the bipotential theory. 
To be short, the space $X = \mathbb{R}^{3}$ is the one of relative velocities 
between points of two bodies, and the space $Y$, identified also to 
$\mathbb{R}^{3}$, is the one of the contact reaction stresses. The duality 
product is the usual scalar  product. We put
$$(x_{n},x_{t})\in X = \mathbb{R} \times \mathbb{R}^{2}, \quad 
    (y_{n},y_{t})\in Y = \mathbb{R} \times \mathbb{R}^{2}\ ,$$
where $x_{n}$ is the gap velocity, $x_{t}$ is the sliding velocity, 
$y_{n}$ is the contact pressure and $y_{t}$ is minus the friction stress. 
The friction coefficient is $\mu > 0$. The graph of the law of unilateral 
contact with Coulomb's dry friction is defined as the union of three sets, 
respectively corresponding to the 'body separation', the 'sticking' and the 
'sliding'. 
\begin{equation}
M = \left\{ (x,0)\in X \times Y \   \mid \ x_{n} < 0 \right\} 
     \cup \left\{ (0,y)\in X \times Y \   \mid \  
                       \parallel y_{t} \parallel \leq \mu y_{n}  \right\} 
		       \cup
\label{Coulomb friction contact law}
\end{equation} 
$$\cup \left\{ (x,y) \in X \times Y \   \mid \  x_{n} = 0, \ x_{t} \neq 0, \ 
                       y_{t} = \mu y_{n} \dfrac{x_{t}}{\parallel x_{t}
		       \parallel}   \right\} $$
It is well known that this graph is not monotone, then not cyclically 
monotone. As usual, we introduce Coulomb's cone
$$ K_{\mu} =  \left\{ (y_{n},y_{t})\in Y \   \mid \  
                      \parallel y_{t} \parallel \leq \mu y_{n}  \right\}  ,$$
and its conjugate cone
$$ K_{\mu}^* =  \left\{ (x_{n},x_{t})\in X \   \mid \  
                      \mu \parallel x_{t} \parallel  + x_{n} \leq 0  \right\}  .$$
In particular, we have
$$ K_{0} =  \left\{ (y_{n},0)\in Y \   \mid \  y_{n} \geq 0  \right\}  , 
\quad 
     K_{0}^* =  \left\{ (x_{n},x_{t})\in X \   \mid \   x_{n} \leq 0  \right\}  .$$
Now, we define some sets useful in the sequel. Let us consider $p > 0$ and the closed 
convex disc obtained by cutting Coulomb's cone at the level $y_{n} = p $
$$ D(p) =  \left\{ y_{t}\in  \mathbb{R}^{2} \   \mid \  
                      \parallel y_{t} \parallel \leq \mu p  \right\}  .$$
Therefore, for each value of $p > 0$, we define a set of 'sticking couples'
$$ M^{\left( a \right) }_{p} =  \left\{ (0,(p,y_{t})) \in X \ \times Y \   \mid \  
                      \  y_{t} \in D (p)  \right\}  \ ,$$
and a set of 'sliding couples'
$$ M^{\left( s \right) }_{p} =  \left\{ ((0,x_{t}),(p,y_{t})) \in X \ \times Y 
\   \mid \  
          \  \parallel y_{t} \parallel = \mu p, \ \exists\lambda > 0, \ 
	  x_{t} = \lambda y_{t}  \right\}  \ .$$
So, we can cover the graph $M$ by the set of following subgraphs parameterized 
by $p \in \left[0,+\infty\right] $
\begin{enumerate}
\item[(a)] $ M_{p} = M^{\left( a \right) }_{p} \cup M^{\left( s \right) }_{p}, 
                                        \quad p\in \left(   0,+\infty \right)  
					\ $,
\item[(b)] $ M_{0} =  \left\{ (x,0)\in X \times Y \   \mid \  x_{n} \leq 0  
\right\}  \ $,
\item[(c)] $ M_{+\infty} = \emptyset , $ by convention.
\end{enumerate}
All these subgraphs are cyclically monotone but none of them is maximal. Let us 
construct by Rockafellar's Theorem the corresponding associated functions  
$\phi_{p}$ and $\psi_{p}$ such that $x_{0}=0$ and 
$\phi_{p}(0) = \psi_{p} (y_{0})=0 $. For $p \in \left( 0,+\infty \right) $, 
the computations give
$$\phi_{p} (x) = p x_{n} + \mu p \parallel x_{t} \parallel \ , 
     \quad \psi_{p} (y) = \chi_{D(p)} (y_{t}) \ .$$
Their Legendre-Fenchel duals are
$$\phi^{*}_{p} (y) = \chi_{ \{ p \} }(y_{n}) + \chi_{D(p)}(y_{t}) \ , 
   \quad \psi ^{*}_{p} (x) = \mu p \parallel x_{t} \parallel + \chi_{ \{ 0 \} } 
   (x_{n}) \ . $$
For $p = 0$, we obtain
$$\phi_{0} (x) = 0 \ , \quad \psi_{0} (y) = \chi_{K_{0}} (y) \ .$$
Their Legendre-Fenchel duals are
$$\phi^{*}_{0} (y) = \chi_{ \{ 0 \} }(y)  \ , \quad \psi ^{*}_{0} (x) =  \chi_{ K^{*}_{0} } (x) \ . $$
For fixed $p$, define the bipotentials $\displaystyle b_{i,p}$, $i=1,2$, by: 
$$b_{1,p}(x,y) = \phi_{p} (x) + \phi^{*}_{p} (y)  \quad , $$
$$b_{2,p}(x,y) = \psi_{p}^{*} (x) + \psi_{p} (y) \quad . $$
As an application of Theorem \ref{thmsup} we obtain that $\displaystyle b_{p} = 
\max \left\{ b_{1,p}, b_{2,p} \right\}$ is a bipotential. Indeed, we shall check
only the point (ii') from Theorem (\ref{thmsup}) (the point (ii") is true by a
similar computation). For $\lambda \in [0,1)$ and $p \not = 0$ we have: 


$$\lambda \phi_{p}(x) + (1-\lambda) \psi_{p}^{*}(x) = \chi_{ \{ 0 \} } (x_{n})
+ \mu p \parallel x_{t} \parallel$$


therefore we get 
$$\left( \lambda \phi_{p}(x) + (1-\lambda) \psi_{p}^{*} \right)^{*}(y) =  
\chi_{D(p)} (y_{t})$$
Also, by computation we obtain: 
$$\lambda \phi_{p}^{*}(y) + (1-\lambda) \psi_{p}(y) = \chi_{ \{ p \} }(y_{n}) +
\chi_{D(p)}(y_{t})$$
If $\displaystyle \phi_{p}^{*}(y) < + \infty \,  , \, \psi_{p}(y) < + \infty$
then in particular $\displaystyle y_{n} = p$ and we obtain (\ref{iiprim}) as an 
equality $0=0$. 
All other cases, involving $\lambda = 1$ or $p = 0$ are solved in the same way. 

The bipotential $b_{p}$ has the expression: 
$$b_{p} (x,y) \ =  \mu p \parallel x_{t} \parallel + \chi_{D(p)}(y_{t}) 
                            + \chi_{ \{ p \} }(y_{n}) + \chi_{ \{ 0 \} } (x_{n}), \quad p\in\left( 0,+\infty \right) \ , $$
$$b_{0} (x,y) \ = \chi_{ \{ 0 \} }(y) + \chi_{ \left( -\infty,0 \right]  } (x_{n}) \ . $$

It is easy to check that the function $\displaystyle p \in [0,+\infty] 
\mapsto b_{p}$ is a bipotential convex cover, therefore by Theorem \ref{thm2} 
we obtain a bipotential for the set $M$. By direct computation, this
bipotential, 
defined as 
$$b(x,y) = \inf \left\{ b_{p}(x,y) \mbox{ : } p \in [0,+\infty] \right\} \quad ,
$$
has the following expression: 
$$b (x,y) \ =  \mu y_{n} \parallel x_{t} \parallel + \chi_{K_{\mu}} (y) +  \chi_{ K^{*}_{0} } (x) \ .$$
Therefore, we recover the bipotential previously given in \cite{saxfeng}. 

\section{Conclusion} 

The present approach shows that  the bipotential related to Coulomb's friction
law  is related to a specific  bipotential  convex cover with the property that 
any graph of the cover is non maximal cyclically monotone.

\vspace{\baselineskip}

\end{document}